\documentstyle[12pt]{article}
\newcommand{\qed} {\hspace {0.1in} \rule {1.5mm} {3.5mm}}

\newtheorem{lemma}{Lemma}[section]
\newtheorem{corollary}{Corollary}[section]
\newtheorem{conjecture}{Conjecture}[section]

\newtheorem{proposition}{Proposition}[section]
\newtheorem{definition}{Definition}[section]
\def\ben{b_E}

\def\1g{1_\Gamma}

\def\dim{{\rm dim}}

\def\<{\langle}
\def\>{\rangle}
\def\proof{\smallskip\noindent{\bf Proof:} }
\def\wk{\widetilde{K}}
\def\wt{\widetilde{T}}
\def\bet{L_{(2)}b}
\def\bF{\mbox{\bf F}_2}
\def\bc{\mbox{\bf C}}
\def\bn{\mbox{\bf N}}
\def\bz{\mbox{\bf Z}}
\def\br{\mbox{\bf R}}
\def\cf{{\it F}}
\def\hep{\frac{\epsilon}{2}}
\def\hg{h_\Gamma}
\def\hta{h^{top}_\alpha}
\def\Fol{F$\mbox{\o}$lner}
\def\l2{\log_2}
\title{Amenable groups, topological entropy and Betti numbers}
\author{{\sc G\'abor Elek}\thanks{Partially supported by OTKA grant
T 25004 and the Bolyai Fellowship}
\cr Alfred Renyi Mathematical Institute of
the Hungarian Academy of Sciences\cr P.O. Box 127, H-1364 Budapest, Hungary}
\date{}
\begin{document}

\maketitle

\noindent{\bf AMS Subject Classifications:} 43A07, 28D20
\vskip 0.2in
\noindent{\bf Keywords:}  amenable groups, topological entropy, group algebras,
$L^2$-Betti numbers
\vskip 0.2in
\noindent{\bf Abstract.} We investigate an analogue of the $L^2$-Betti numbers
for amenable linear subshifts. The role of the von Neumann dimension
shall be played by the topological entropy.
\newpage
\section{Introduction}
Let $\Gamma$ be a finitely generated group. Then the Hilbert space
$l_2(\Gamma)$ has a natural left $\Gamma$-action by translations:
$$L_\gamma(f)(\delta)=f(\gamma^{-1}\delta)\,. $$
Using the so-called von Neumann dimension we can assign a real number to any
$\Gamma$-invariant linear subspace of $[l^2(\Gamma)]^n, n\in \bn$ satisfying
the following basic axioms \cite{Pan}.

\begin{enumerate}
\item {\bf Positivity:} If $V\subset [l^2(\Gamma)]^n$ $\Gamma$-invariant
linear subspace, then $\dim_\Gamma (V)\geq 0$. Also, $\dim_\Gamma (V)=0$
if and only if $V=0$.
\vskip 0.2in
\item {\bf Invariance:} If $V\subset [l^2(\Gamma)]^n$,
 $W\subset [l^2(\Gamma)]^m$ and $T$ is a $\Gamma$-equivariant isomorphism
from $V$ to a dense subset of $W$, then $\dim_\Gamma (V)= \dim_\Gamma (W)$.
\vskip 0.2in
\item {\bf Additivity:} If $Z$ is the orthogonal direct sum of $V$ and $W$,
 then $\dim_\Gamma (Z)=\dim_\Gamma (V)+\dim_\Gamma (W)$.
\vskip 0.2in
\item {\bf Continuity:} If $V_1\supset V_2\supset \dots$ is
a decreasing sequence of $\Gamma$-invariant linear subspaces, then:
$$\dim_\Gamma(\cap_{j=1}^\infty V_j)=\lim_{j\rightarrow \infty}
\dim_\Gamma (V_j) \,.$$
\vskip 0.2in
\item {\bf Normalization:} $\dim_\Gamma [l^2(\Gamma)]=1$.
\vskip 0.2in

\end{enumerate}
There is an important application of the von Neumann dimension in algebraic
topology due to Atiyah \cite{At} (see also \cite{Dodz}). He defined certain 
invariants of finite simplicial complexes : the $L^2$-Betti numbers. The idea
is the following, let $\widetilde{K}$ be an infinite, simplicial complex with
 a free and simplicial $\Gamma$-action as covering transformations such that
$\widetilde{K}/\Gamma=K$ is finite. Denote by $C^p_{(2)}(\wk)$ the Hilbert
space of square-summable, real $p$-cochains of $\wk$.
 Then one has the following
differential complex of Hilbert spaces,
$$ C^0_{(2)}(\wk)\stackrel{d_0}{\rightarrow}
C^1_{(2)}(\wk)\stackrel{d_1}{\rightarrow}\dots\stackrel{d_{n-1}}{\rightarrow}
C^n_{(2)}(\wk)\,,$$
where the $d_p$'s are the usual coboundary operators. Note that
$C^p_{(2)}(\wk)\cong [l^2(\Gamma)]^{|K_p|}$, where $K_p$ denotes the set of
$p$-simplices in $K$. Atiyah's $L^2$-Betti numbers are defined as
$$\bet^p(K)=\dim_\Gamma \,Ker\,d_p-\dim_\Gamma\,Im\,d_{p-1}\,.$$
Let us list some basic results on the $L^2$-Betti numbers.
\begin{itemize}
\item (Dodziuk, \cite{Dodz}) If $\wk$ and $\widetilde{L}$ are
homotopic by a $\Gamma$-invariant homotopy, then the
corresponding $L^2$-Betti numbers of $\wk/\Gamma=K$ and
 $\widetilde{L}/\Gamma=L$ are equal.
\item (Cohen, \cite{Coh})
$$\sum_{p=0}^n (-1)^p\bet^p(K)=e(K)\,,$$
the Euler characteristic of $K$.
\item (Cheeger \& Gromov, \cite{CG})
If $\wk$ is contractible and $\Gamma$ is amenable, then
all $L^2$-Betti numbers are vanishing.
\item (Linnell,  \cite{Lin}) If $\Gamma$ is elementary amenable and torsion-free
then all $L^2$-Betti numbers are integers.
\item (L\~uck, \cite{Luck2}) Let $\Gamma$ be residually finite and
$$\Gamma\supset\Gamma_1\supset\Gamma_2\dots,
\quad\cap^\infty_{i=1}\Gamma_i=1_\Gamma$$
normal subgroups of finite index and let $X_i=\wk/\Gamma_i$ the
corresponding finite coverings of $K$. Then
$$\bet^p(K)=\lim_{i\rightarrow\infty} 
\frac {\dim_{\br} H^p(X_i,\br)}{|\Gamma:\Gamma_i|}$$
\item (Dodziuk \& Mathai \cite{DM} ) If $\{L_n\}_{n=1}^\infty$
 is an exhaustion of $\wk$
by finite simplicial complexes spanned by a $\{\cf_n\}_{n=1}^\infty$
\Fol-exhaustion  , then
$$\bet^p(K)=\lim_{i\rightarrow\infty} 
\frac {\dim_{\br} H^p(L_n,\br)}{|\cf_n|}$$
\end{itemize}
Note that the second and the third results together imply that if
$K$ is an acyclic simplicial complex with amenable fundamental group
then its Euler characteristic is zero \cite{CG}. Another interesting
application is due to L\"uck: If $\Gamma$ is amenable, then
the group algebra $\bc[\Gamma]$ as a free module over itself
generates an infinite cyclic subgroup in the Grothendieck group
of $\bc[\Gamma]$ \cite{Luck}.
\\
\\
The analogue setting we are investigating in this paper is the following.
Let $\Gamma$ be a finitely generated  amenable group
(see \cite{EG} why amenability is crucial). We denote by
 $\sum_\Gamma$ the full Bernoulli shift that is the linear space of
$\bF$-valued functions on $\Gamma$, where $\bF$ is the field of two elements.
The space $\sum_\Gamma$ is a compact, metrizable space in the pointwise
 convergence topology
equipped with the natural left $\Gamma$-action by translations. A space
$V\subset[\sum_\Gamma]^n$ is a linear subshift if it is linear as a
$\bF$-vector space, closed in the topology and invariant with
respect to the $\Gamma$-action. The notion of dimension
is the topological entropy of the linear subshifts. This is well-known for
$\bz$ and $\bz^d$-actions and somehow less-known for general amenable group
actions (nevertheless see \cite{MO}).
 We shall observe that our dimension
$h_\Gamma$ satisfies similar axioms as $\dim_\Gamma$ :
\begin{enumerate}
\item {\bf Nonnegativity:} For any $V$ linear subshift : $\hg(V)\geq 0$. But it can be
 zero even if $V$ is not zero.
\vskip 0.2in
\item{\bf Monotonicity:}  If $V\subset W$, then $\hg(V)\leq \hg(W)$.
\vskip 0.2in
\item{\bf Invariance:}  If $T:V\rightarrow W$ continuous $\Gamma$-equivariant
 linear
isomorphism, then $\hg (V)=\hg (W)$.
\vskip 0.2in
\item{\bf Additivity:} If $ Z=V\oplus W$ , then $\hg(Z)=\hg(V)+\hg(W)\,.$
\vskip 0.2in
\item{\bf Continuity:} If $V_1\supset V_2\dots $ is a decreasing sequence of linear
subshifts, then :
$$\hg(\cap_{j=1}^\infty V_j)=\lim_{j\rightarrow \infty}
h_\Gamma (V_j) \,.$$
\item {\bf Normalization:} $\hg(\sum_\Gamma)=1$.
\end{enumerate}
Now let $\wk$ be as above. Then
we have the ordinary cochain complex of
$\bF$-coefficients over $\wk$:
$$ C^0(\wk,\bF)\stackrel{d_0}{\rightarrow}
C^1(\wk,\bF)\stackrel{d_1}{\rightarrow}\dots\stackrel{d_{n-1}}{\rightarrow}
C^n(\wk,\bF)\,.$$
Then the $p$-cochain space $C^p(\wk,\bF)$ is $\Gamma$-isomorphic
to $[\sum_\Gamma]^{|K_p|}$, where $K_p$ denotes
the set of $p$-simplices in $K$.
We define the $p$-th entropy Betti number $\ben^p(K)$ as
$\hg(Ker\,d_p)-\hg(Im\,d_{p-1})$.
In this paper we shall prove the following analogues of the
$L^2$-results.
\begin{itemize}
\item  If $\wk$ and $\widetilde{L}$ are
homotopic by a $\Gamma$-invariant homotopy and $\Gamma$
is poly-cyclic then the
corresponding entropy-Betti numbers of $\wk/\Gamma=K$ and
 $\widetilde{L}/\Gamma=L$ are equal.
\item 
$$\sum_{p=0}^n (-1)^p\ben^p(K)=e(K)\,,$$
the Euler characteristic of $K$.
\item 
If $\wk$ is contractible  then
all entropy-Betti numbers are vanishing.
(this is quite obvious, the point is that the corollary on
the vanishing Euler-characteristic still follows from this and the
previous statement)
\item  If $\Gamma$ is poly-infinite-cyclic
then all entropy-Betti numbers are integers.
\item  Let $\Gamma$ be free Abelian and
$$\Gamma\supset\Gamma_1\supset\Gamma_2\dots,
\quad\cap^\infty_{i=1}\Gamma_i=1_\Gamma$$
normal subgroups of finite index and let $X_i=\wk/\Gamma_i$ the
corresponding finite coverings of $K$. Then
$$\ben^p(K)=\lim_{i\rightarrow\infty} 
\frac {\dim_{\bF} H^p(X_i,\bF)}{|\Gamma:\Gamma_i|}$$
\item  If $\{L_n\}_{n=1}^\infty$ is an exhaustion of $\wk$
by finite simplicial complexes spanned by a $\{\cf_n\}_{n=1}^\infty$ \Fol-exhaustion, then
$$\ben^p(K)=\lim_{n\rightarrow\infty} 
\frac {\dim_{\bF} H^p(L_n,\bF)}{|\cf_n|}$$
\end{itemize}
We shall also prove an analogue of L\"uck's result on the Grothendieck-group
for the group algebras $\bF[\Gamma]$.
\section{Amenable groups and quasi-tiles}
Let $\Gamma$ be a finitely generated group with a symmetric generator
set $\{g_1,g_2,\dots,g_k\}$. The right Cayley-graph of $\Gamma$, $C_\Gamma$ is
defined as follows. Let $V(C_\Gamma)=\Gamma$, $E(C_\Gamma)=\{(a,b)\in
\Gamma\times\Gamma:\,\, $ there exists $g_i:\, ag_i=b\}$. The shortest
path distance $d$ of $C_\Gamma$ makes $\Gamma$ a discrete metric space.
We shall use the following notation. If $H\subset\Gamma$ is a finite set, then
$B_r(H)$ is the set of elements $a$ in $\Gamma$ such that there exists
$h\in H, d(a,h)\leq r$. We denote $B_1(H)\backslash H$ by $\partial H$ and
$B_r(H)\backslash H$ by $\partial_r H$. An exhaustion of $\Gamma$ by finite 
sets 
$$1_\Gamma\in
\cf_1\subset \cf_2 \subset \dots ,\quad \cup_{j=1}^\infty \cf_j=\Gamma$$
is called a \Fol-exhaustion if for any $r\in \bn$:
$\lim_{n\rightarrow\infty}\frac {|\partial_r \cf_n|} {|\cf_n|}=0$.
A group $\Gamma$ is called amenable if it possess a \Fol-exhaustion.
Some amenable groups have {\it tiling} \Fol-exhaustion that is any $\cf_n$
is a tile : There exists $C\subset\Gamma$ such that $\{c\cf_n\}_{c\in C}$
is a partition of $\Gamma$. For example $\bz^n$ has this tiling property.
As observed by Ornstein and Weiss \cite{OW} {\it any} amenable
 group has quasi-tiling
\Fol-exhaustion. Let us recall their construction. Let $\{A_i\}_{i=1}^\infty$
be finite sets. Then we call them $\epsilon$-disjoint if there exist subsets
$\overline{A_i}\subset A_i$ so that $\overline{A_i}\cap\overline{A_j}=0$
if $i\neq j$, and $\frac {|\overline{A_i}|}{|A_i|}\geq 1-\epsilon$ for all
$i$. Now let $B$ another finite set. We say that $\{A_i\}_{i=1}^\infty$
$(1-\epsilon)$-cover $B$, if
$$\frac {|B\cap \cup_{i=1}^\infty A_i|}{|B|}\geq  1-\epsilon\,.$$
The subsets of $\Gamma$, $\1g\in T_1\subset T_2\subset\dots \subset T_N$
form an $\epsilon$-quasi-tile system if for any finite subset of
$A\subset\Gamma$, there exists $C_i\subset\Gamma$, $i=1,2,\dots,N$ such that

1. $C_iT_i\cap C_jT_j=\emptyset$ if $i\neq j$.

2. $ \{cT_i:\, c\in C_i\}$ are $\epsilon$-disjoint sets for any fixed $i$.

3. $\{C_iT_i\}$ form a $(1-\epsilon)$-cover of $A$.

\noindent
The following proposition is Theorem 6. in \cite{OW}:
\begin{proposition}\label{ow}
If $\cf_1\subset\cf_2\subset\dots $ is a \Fol-exhaustion of an amenable group, then
for any $\epsilon>0$ we can choose a finite subset
$\cf_{n_1}\subset \cf_{n_2}\subset
\dots\subset\cf_{n_N}$ such that they form an
$\epsilon$-quasi-tile system. The number $N$ may depend on $\epsilon$.
\end{proposition}

\section{The topological entropy of linear subshifts}
First of all we define an averaged dimension $\hg(W)$ for
linear subshifts and then we shall show that it coincides with the
topological entropy.
Let $\Gamma$ be a finitely generated amenable group
with \Fol-exhaustion $1_\Gamma\in\cf_1\subset \cf_2 \subset \dots ,\quad
\cup_{j=1}^\infty \cf_j=\Gamma$. We introduce some notations.
If $\Lambda\subset\Gamma$ is a finite set, then let
$[\sum_\Lambda]^r$ be the space of functions in $[\sum_\Gamma]^r$
supported on $\Lambda$. Also, $[\sum_\Gamma^0]^r$ denotes the space
 of finitely supported functions.
Now let $W\subset [\sum_\Gamma]^r$ be a $\Gamma$-invariant not
 necessarily closed linear subspace. Then for any finite
$\Lambda\subset\Gamma$
let $W_\Lambda\subset [\sum_\Lambda]^r$ be the linear space of
 functions $\eta$ supported on $\Lambda$ such that there exists
$\nu\in W$ : $\eta\mid_\Lambda=\nu\mid_\Lambda$.
\begin{definition}
$\hg(W)=\limsup_{n\rightarrow\infty}\frac {\l2|W_{\cf_n}|}{|\cf_n|}$
\end{definition}
Note that $\l2|W_{\cf_n}|$ is just the dimension of the vector space
$W_{\cf_n}$ over the field $\bF$. It will be obvious from the next proposition
that $\hg(W)$ does not depend on the particular choice of the exhaustion.
\begin{proposition} 
{}

\begin{enumerate}\item
  $\hg(W)=\hg(\overline{W})$, where $\overline{W}$ denotes
the closure of $W$ in the pointwise convergence topology.
\item
  $\hg(W)=\liminf_{n\rightarrow\infty}\frac {\l2|W_{\cf_n}|}{|\cf_n|}$,
hence $\lim_{n\rightarrow\infty}\frac {\l2|W_{\cf_n}|}{|\cf_n|}$
always exists and equals to $\hg(W)$.\end{enumerate}
\end{proposition}
\proof The first part is obvious from the definition, for the second part
we argue by contradiction. Suppose that
$$\hg(W)-\liminf_{n\rightarrow\infty}\frac
{\l2|W_{\cf_n}|}{|\cf_n|}=\delta>
0.$$
Consider a subsequence $\cf_{n_1}\subset \cf_{n_2}\subset\dots $ such that $$
\sup_{i\rightarrow\infty}\frac {\l2|W_{\cf_{n_i}}|}{|\cf_{n_i}|}\leq
\liminf_{n\rightarrow\infty}\frac {\l2|W_{\cf_n}|}{|\cf_n|}+\epsilon\,,$$
where the explicite value of $\epsilon$ shall be chosen later accordingly.
Then pick an $\epsilon$-quasi-tile system from our subsequence:
$\cf_{m_1}\subset \cf_{m_2}\subset\dots \subset \cf_{m_N}$. Now we take
 an arbitrary $\cf_n$ from the original \Fol-exhaustion. By Proposition
\ref{ow} we have an $\epsilon$-disjoint, $(1-\epsilon)$-covering
of $\cf_n$ by translates of the quasi-tile system. Denote by $R_1, R_2,\dots
R_k$ those tiles which are properly contained in $\cf_n$. Then we have 
the following estimate.
\begin{equation}
\label{121}
|W_{\cf_n}|\leq 2^{(r\epsilon|\cf_n|+r|\cf_n\backslash
 B_{D+1}(\partial\cf_n)|)}\prod_{i=1}^k|W_{R_i}|,
\end{equation}
where $D$ is the diameter of the largest tile $\cf_{m_N}$. The inequality
(\ref{121}) follows from the fact that a function $\xi\in W_{\cf_n}$ is
uniquely determined by its restrictions on the covering tiles and its
restriction on the uncovered elements. The later one consists of two
parts; the elements which are not covered by the original covering and
the elements which are covered by tiles intersecting the complement of
$\cf_n$.
These ``badly'' covered elements are in a $D+1$-neighbourhood of the boundary
of $\cf_n$. Also, by $\epsilon$-disjointness we have the estimate
\begin{equation}
\label{131}
\sum_{i=1}^k|R_i|\leq\frac{1}{1-\epsilon} |\cf_n|\,.
\end{equation}
Therefore,
$$|W_{\cf_n}|\leq 2^{(r\epsilon|\cf_n|+r|\cf_n\backslash
 B_{D+1}(\partial\cf_n)|)}2^{\frac{1}{1-\epsilon}
(\hg-\delta+\epsilon)|\cf_n|}$$
Hence,
\begin{equation}
\frac{\log_2\,|W_{\cf_n}|}{\cf_n}\leq 
r\epsilon+\frac{r|\cf_n\backslash
 B_{D+1}(\partial\cf_n)|)}{|\cf_n|}+\frac{1}{1-\epsilon}(\hg-\delta+\epsilon)
\end{equation}
Consequently if we choose $\epsilon$ small enough, then for large $n$,
$\frac{\log_2\,|W_{\cf_n}|}{|\cf_n|}\leq\hg-\frac{\delta}{2}$, leading
 to a contradiction.\qed

Now we recall the notion of topological entropy. Let $\Gamma$ be an amenable
group as above
and let $X$ be a compact metric space equipped with a continuous
$\Gamma$-action; $\alpha:\Gamma\rightarrow\,Homeo(X)$. Instead of
 the original definition of Moulin-Ollagnier \cite{MO} we use the
equivalent ``spanning-separating'' definition, that is a direct
generalization of the Abelian case \cite{KS}. We call a finite set $S\subset X$
$(n,\epsilon)$-separated if for any distinct points $s,t\in S$ there exists
$\gamma\in \cf_n$ such that $d(\alpha(\gamma^{-1})(s),\alpha(\gamma^{-1})(t))>
\epsilon$. We denote by $s(n,\epsilon)$ the maximal
cardinality of such sets.
We call a finite set $R
\subset X$ $(n,\epsilon)$-spanning if for any $x\in X$ there exists
$y\in R$ such that $d(\alpha(\gamma^{-1})(x),\alpha(\gamma^{-1})(y))\leq
\epsilon$, for all $\gamma\in \cf_n$.  We denote by $r(n,\epsilon)$ the minimal
cardinality of such sets. Obviously if $\epsilon'<\epsilon$ then 
$s(n\epsilon')\geq s(n\epsilon), r(n\epsilon')\geq r(n\epsilon)$. Also,
we have the inequalities:
$$r(n,\epsilon)\leq s(n,\epsilon)\leq r(n,\hep)\,.$$
Indeed, any $(n,\epsilon)$-separating set is $(n,\epsilon)$-spanning.
On the other hand if $R=x_1,x_2,\dots,x_k$ is a $(n,\hep)$-spanning
set then $X=\bigcup_{i=1}^k\, D(x_i,n,\hep)$, where
$$D(x_i,n,\hep)=\{y\in
X:\,\,d(\alpha(\gamma^{-1})(x),\alpha(\gamma^{-1})(y))\leq\hep,
\mbox{for all $\gamma\in \cf_n$}\}\,.$$
Any $D(x_i,n,\hep)$ can contain at most one element of
a $(n,\epsilon)$-separating set, hence $s(n,\epsilon)\leq r(n,\hep)$.

\noindent
Consequently, $\lim_{\epsilon\rightarrow 0} \limsup_{n\rightarrow\infty} \l2
r(n,\epsilon)=\lim_{\epsilon\rightarrow 0} \limsup_{n\rightarrow\infty} \l2
s(n,\epsilon)$. This joint limit is called the topological entropy of
the $\Gamma$-action and denoted by $h^{top}_\alpha(X)$.
Note that it follows from the definition that $\hta(X)$ depends
only on the topology and not the particular choice of the metric on $X$.
\begin{proposition}
\label{331}
Let $V\subset[\sum_\Gamma]^r$ be a linear subshift. Then
$h^{top}_L(V)=\hg(V)$. 
\end{proposition}
\proof First we fix a metric on $V$ that defines the pointwise convergence
topology.
If $v,w\in V$, then let $d_V(v,w)=2^{-(n-1)}$, where $n$ is the infinum of $k$'s
such that $v\mid_{\cf_k}\neq w\mid_{\cf_k}$.
First note that $|V_{\cf_n}|\leq s(n,1)$. Indeed
 if $v_1,v_2,\dots v_s$ is a subset of $V$ such that
$v_i\mid_{\cf_n}\neq v_j\mid_{\cf_n}$ when $i\neq j$, then
there exists $\gamma\in\cf_n$ such that
$L_{\gamma^{-1}}v_i(1_\Gamma)\neq L_{\gamma^{-1}}v_j(1_\Gamma)$.
 Fix an $\epsilon$
and choose $K_\epsilon,M_\epsilon\in \bn$ such that
$\epsilon>2^{-K_\epsilon}$ and $\cf_{K_\epsilon}\subset
B_{M\epsilon}(1_\Gamma)$. Then we claim that $s(n,\epsilon)\leq
|V_{B_{M\epsilon}(\cf_n)}|$. Indeed, if $x\mid_{B_{M\epsilon}(\cf_n)}=
y\mid_{B_{M\epsilon}(\cf_n)}$, then for any $\gamma\in\cf_n$,
$d_V(L_{\gamma^{-1}}(x),L_{\gamma^{-1}}(y))\leq 2^{-K_\epsilon}<\epsilon$.
Hence if $\epsilon<1$,
$$\hg(V)=\lim_{n\rightarrow\infty}\frac{\l2|V_{\cf_n}|}{|\cf_n|}
\leq \l2(s(n,\epsilon))\leq \lim_{n\rightarrow\infty}
\frac{\l2|V_{B_{M\epsilon}(\cf_n)}|}{|\cf_n|}=\hg(V)\,\,\,\qed$$

Note that the previous proposition immediately shows that $\hg(V)$ does not
depend on how $V$ is imbedded in a full shift.
\section{Extended Configurations}
The notion of extended configuration is due to Ruelle in a slightly
different form. Again we start with
a linear subshift $V\subset[\sum_\gamma]^r$. For
any $\Lambda\subset \Gamma$ finite set
 let $V^\Omega_\Lambda\subset[\sum_\gamma]^r$ be a finite dimensional
 linear subspace satisfying the following
axioms:
\begin{itemize}
\item {\bf Extension:} $V_\Lambda\subset V^\Omega_\Lambda$.
\item {\bf Invariance:} $V^\Omega_{\gamma\Lambda}=L_\gamma(V^\Omega_\Lambda)$.
\item {\bf Transitivity:} If $\Lambda\subset M$, then for any $\xi\in
V^\Omega_M$ there exists $\mu\in V^\Omega_\Lambda$ such that
$\xi\mid_\Lambda=\mu\mid_\Lambda$.
\item {\bf Determination:} If $\xi\in [\sum_\gamma]^r$ and for any
$\Lambda\subset\Gamma$ finite there exists $\xi_\Lambda\in V^\Omega_\Lambda$
such that $\xi\mid_\Lambda=\xi_\Lambda$. \end{itemize}
We call such a system an extended configuration of $V$. Its topological entropy
 is defined as $\hg^\Omega(V)=\limsup_{n\rightarrow\infty}\frac
{|V^\Omega_{\cf_n}|}{|\cf_n|}\,.$
\begin{proposition}
\label{151}
$h^\Omega_\Gamma(V)=\hg(V)$ 
(compare to Theorem 3.6 \cite{Rue})
\end{proposition}
\proof Let us suppose that $h^\Omega_\Gamma(V)-\hg(V)=\delta>0$.
 Again choose an $\epsilon$-quasi-tile system $\cf_{n_1},\cf_{n_2}\dots,
\cf_{n_N}$ such that
$|V_{\cf_{n_i}}|<2^{(\hg +\epsilon)|\cf_{n_i}|}$, where the explicite
 value of $\epsilon$ will be given later. Denote by $V^k_{\cf_{n_i}}$
the space of functions $\mu$ in $V_{\cf_{n_i}}$ such that there exists
$\xi\in V^\Omega_{B_k(\cf_{n_i})}$ with $\xi\mid_{ \cf_{n_i}}=\mu\mid_
 { \cf_{n_i}}$. By Extension and Determination properties it is easy to see
 that for large $p$ :
\begin{equation}
\label{152}
V^p_{\cf_{n_i}}=V_{\cf_{n_i}},
\end{equation}
where $1\leq i \leq N$.
Let us pick  a large $p$. Then we can proceed almost the same way as in 
the previous section. Denote by $R_i$, $1\leq i\leq k$ those translates
in the $\epsilon$-disjoint $(1-\epsilon)$-covering of a \Fol-set $\cf_n$ such 
that not only the $R_i$'s but even the $B_p(R_i)$ balls are  contained
in $\cf_n$.
Then by Transitivity we have the following estimate
$$|V^\Omega_{\cf_n}|\leq
2^{(r\epsilon |{\cf_n}|+r|\cf_n\backslash B_{p+D+1}(\partial\cf_n)|)}
\prod_{i=1}^k|V_{R_i}^p| \quad.$$
That is by Invariance, (\ref{131}) and (\ref{152}) :
$$\frac{\l2 |V^\Omega_{\cf_n}|}{|\cf_n|}\leq
r\epsilon +
\frac{r|\cf_n\backslash B_{p+D+1}(\partial\cf_n)|}{|{\cf_n}|}
+\frac{1}{1-\epsilon}(\hg^\Omega-\delta+\epsilon)$$
which leads to a contradiction provided that we choose $\epsilon$ small
enough.\qed
\section{Basic Properties}
Now we are in the position to prove the basic properties of $\hg$ as stated in 
the Introduction. The Monotonicity, Normalization and Positivity axioms
are obviously satisfied.
\begin{lemma}
\label{161}
Let $V,W\subset[\sum_\Gamma]^r$ be linear subshifts such that $V\cap W=0$,
 then $\hg(V)+\hg(W)=\hg(V\oplus W)\,.$
\end{lemma}
\proof
First note that just because $V\cap W$ is the zero subspace
it is not necessarily true that $V_\Lambda\cap W_\Lambda=0$ as well.
However, we can prove that $N^\Omega_\Lambda=V_\Lambda\cap W_\Lambda$
is an extended configutation of the zero subspace. We only need to
 show that the Determination axiom is satisfied. 
Suppose that $\xi\in [\sum_\Gamma]^r$ such that
$\xi\mid_{\cf_n}\in V_{\cf_n}\cap W_{\cf_n}$. Therefore there exists
$v_n\in V, w_n\in W$  such that $\xi\mid_{\cf_n}=v_n\mid_{\cf_n}=
w_n\mid_{\cf_n}$. Hence $v_n\rightarrow \xi, w_n\rightarrow \xi$
in the topology of $\sum_\Gamma)^r$. The spaces $V$ and $W$ are
closed, thus $\xi\in V\cap W$, hence $\xi=0$. By elementary linear
algebra,
$$\dim_{\bF}(V\oplus W)_\Lambda + \dim_{\bF} N^\Omega_\Lambda=
\dim_{\bF}V_\Lambda + \dim_{\bF}W_\Lambda\quad.$$
Hence by our Proposition \ref{151} our Lemma follows.\qed

\noindent
Now we prove a property of the entropy that is slightly more
general then invariance.
\begin{proposition}
\label{181}
Let $T:V\rightarrow W$ be a continuous $\Gamma$-equivariant
linear map between linear subshifts $V\subset [\sum_\Gamma]^r, 
W\subset [\sum_\Gamma]^r$. Then $\hg(Ker\, T)+\hg(Im\,T)=\hg(V)$.
\end{proposition}
\proof
First of all let us note that by the compactness of $V$ and the continuity of
$T$ both $Ker\, T$ and $Im\, T$ are linear subshifts. Now let 
us consider the natural right action of $\bF(\Gamma)$ on $\sum_\Gamma$,
$R_\gamma f(x)=f(x\gamma)$. This action obviously commutes 
with our previously defined left $\Gamma$-action. The right action can be
extended
to $s\times r$-matrices with coefficients in $\bF(\Gamma)$ acting
on the column vectors $[\sum_\Gamma]^r$. Obviously, any such
matrix $M$ defines a $\Gamma$-equivariant map; $T_M:[\sum_\Gamma]^r
\rightarrow [\sum_\Gamma]^s$.
\begin{lemma}
\label{191}
Any continuous $\Gamma$-equivariant
linear map $T:V\rightarrow W$ can be given
via multiplication by some $s\times r$-matrix $T_M$ with
coefficients in $\bF(\Gamma)$.
\end{lemma}
\proof:
Since $T$ is uniformly continuous
the value of $T(v)(1_\Gamma)$ is determined by the
value of $v$ on a finite ball $B$, where $B$ does not depend
on $v$. Hence for
any $1\leq i \leq s$ :
$$T(v)(\1g)=\sum_{\gamma\in B} \sum_{j=1}^r c_{ij}^\gamma\cdot
v_j(\gamma)\quad,$$
where $c_{ij}^\gamma\in\bF$. By  $\{T_M\}_{ij}=
\{\sum_{\gamma\in B}c_{ij}^\gamma\gamma\}\in Mat_{s\times r}(\bF[\Gamma])$
define a $s\times r$-matrix. Then for any $v\in V, T_M(v)(\1g)=T(v)(\1g)$.
Hence by the $\Gamma$-equivariance of the matrix multiplication:
$$T_M(v)(\gamma)= L_{\gamma^{-1}}(T_M(v))(\1g)= T_M(L_{\gamma^{-1}}(v))(\1g)
=T(L_{\gamma^{-1}}(v))(\1g)=T(v)(\gamma)\quad\qed$$

\noindent
Now we return to the proof of our Proposition. We denote by $T_M$ the
matrix and by $k$ the diameter of the ball $B$ defined in our Lemma.
Let
$$N^\Omega_\Lambda=\{\mbox{$ v\in [\sum_\Lambda]^r$}: \mbox { there exists}\,
z\in V_{B_k(\Lambda)}, \mbox { such that}\, z\mid_\Lambda=v \mbox { and}\,
T(z)\mid_\Lambda=0\} $$
$$M^\Omega_\Lambda=\{\mbox{$ w\in [\sum_\Lambda]^s$}: \mbox { there exists}\,
z\in V_{B_k(\Lambda)}, \mbox { such that}\, T(z)\mid_\Lambda=w\}$$
Then $N^\Omega_\Lambda$ is an extended configuration of
$Ker\,T_M$, $M^\Omega_\Lambda$ is an extended configuration
of $Im\,T_M$. Let $\wt_\Lambda:
V_{B_k(\Lambda)}\rightarrow [\sum_\Lambda]^s$
be the restriction of $T$ onto $\Lambda$. Then $Im\,\wt_\Lambda=
M^\Omega_\Lambda$. We have the usual pigeon-hole estimate :
\begin{equation}
\label{201}
|N^\Omega_{\cf_n}|\leq |Ker\, \wt_{\cf_n}|\leq |N^\Omega_{\cf_n}|
2^{r|B_k(\partial\cf_n)|}
\end{equation}
Also, by linear algebra we obtain
$$\dim_{\bF} Ker\, \wt_{\cf_n}+ \dim_{\bF} Im\, \wt_{\cf_n}=
\dim_{\bF} V_{B_k(\cf_n)}$$
That is
\begin{equation}
\label{202}
\l2|Ker\, \wt_{\cf_n}|+\l2|Im\, \wt_{\cf_n}|=\l2 |V_{B_k(\cf_n)}|
\end{equation}
It is easy to see that (\ref{201}) and (\ref{202}) together
implies the statement of our Proposition. \qed

\noindent
Now we prove the Continuity property.
\begin{proposition}
\label{211}
If $V^1\supset V^2\supset\dots$ is a decreasing sequence
of linear subshifts then
$$\hg(\bigcap_{j=1}^\infty V^j)=\lim_{j\rightarrow\infty} \hg(V^j)\quad.$$
\end{proposition}
\proof
For any $k\in N, V^k\supset \bigcap_{j=1}^\infty V^j$. Hence
$\hg(\bigcap_{j=1}^\infty V^j)\leq\lim_{j\rightarrow\infty} \hg(V^j)$.
We need to prove the converse inequality. Suppose that for all $j$,
$\hg(V^j)\geq \hg(\bigcap_{j=1}^\infty V^j)+2\epsilon$. Let $T(j)$ be a monotone
increasing function such that
$$\frac {\l2|V^j_{\cf_s}|}{|\cf_s|}\geq
\hg(\bigcap_{j=1}^\infty V^j)+\epsilon\quad,$$
when $s\geq T(j)$. We define a monotone non-increasing function
$S:\bn\rightarrow\bn$ such that $S(n)=1$ if there is no such $j$ so that
$T(j)\leq n$ and $S(n)=\inf\,\{j: T(j)\leq n\}$ otherwise. Then
$S(n)\rightarrow\infty$ as $n\rightarrow\infty$.
Define $V_\Lambda^\Omega=V^{S(n)}_\Lambda$, where $n$
is the smallest integer such that
$|\cf_n|\geq |\Lambda|$. It is easy to see that $V_\Lambda^\Omega$
is an extended configuration of
$\bigcap_{j=1}^\infty V^j$. Therefore,
$\frac {\l2|V^\Omega_{\cf_n}|}{|\cf_n|}\rightarrow\hg(\bigcap_{j=1}^\infty V^j)$.
On the other hand, by our construction: $\frac {\l2|V^\Omega_{\cf_n}|}{|\cf_n|}
\geq \hg(\bigcap_{j=1}^\infty V^j)+\epsilon$, leading
 to a contradiction.\qed
\section{Pontryagin Duality}
In this section we recall the Pontryagin Duality theory \cite{HR}.
Let $A$ be a locally compact Abelian group and let $\widehat{A}$
be its dual. That is the group of continuous homomorphisms
$\chi:A\rightarrow S^1=\{z\in C:\,\,|z|=1\}$. According to the duality theorem
$A$ is naturally isomorphic to its double dual. The relevant example
for us is $A=[\sum_\Gamma]^r$, its dual is $[\sum_\Gamma^0]^r$ the
group of finitely supported elements, where $\langle\chi,f\rangle=
\sum_{\gamma\in\Gamma}(\chi(\gamma),f(\gamma))$ for $\chi\in[\sum_\Gamma^0]^r$
and $f\in [\sum_\Gamma]^r$. Here $(\underline{a},\underline{b})$ is defined
as $\sum_{i=1}^r a_ib_i$. The additive group of $\bF$ is viewed as the
subgroup
$\{-1,1\}\in S^1$. If $H\subset [\sum_\Gamma]^r$ is a compact subgroup
then
$$H^\perp=\{\mbox{$\chi\in [\sum_\Gamma^0]^r:$}\,\langle\chi,h\rangle=1 \mbox
{for any $h\in H$}\}$$
Conversely if $B\subset [\sum_\Gamma^0]^r$ is a subgroup
then
$B^\perp=\{f\in [\sum_\Gamma]^r:\,\langle\chi,f\rangle=1\,\, \mbox
{for any\, $\chi\in B$}\}\quad.$
Then $(H^\perp)^\perp=H, \, (B^\perp)^\perp=B$. If $A,B$ locally compact
groups and $\psi:A\rightarrow B$ continuous homomorpisms, then
its dual $\widehat{\psi}:\widehat{B}\rightarrow
\widehat{A}$ is defined by $\langle\widehat{\psi}(\chi),a\rangle=
\langle\chi,\psi(a)\rangle$. Again the double dual of $\psi$ is itself
if $A$ and $B$ are both compact or both discrete. Then $\psi$ is injective
resp. surjective if and only if $\widehat{\psi}$ is surjective
(resp. injective). Moreover, if we have a short exact sequence
of compact or discrete groups
$$1\rightarrow A_1\rightarrow A_2\rightarrow\dots\rightarrow A_n\rightarrow
 1$$
Then its dual sequence
$$0\rightarrow\widehat{A}_n\rightarrow\dots\widehat{A}_2\rightarrow
\widehat{A}_1\rightarrow 0$$
is also exact.
The next proposition is a version of a result of Schmidt \cite{KS}.
\begin{proposition}
\label{361}
The Pontryagin duality provides a one-to-one correspondance between linear
subshifts and finitely generated left $\bF[\Gamma]$-modules.
\end{proposition}
\proof
First note that if $L_\gamma$ is the left multiplication
by $\gamma$ on $[\sum_\Gamma]^r$ then $\widehat{L}_\gamma$ is the
left multiplication by $\gamma^{-1}$ on $[\sum_\Gamma^0]^r$.
Hence if $V\stackrel{i}{\rightarrow} [\sum_\Gamma]^r$ is
 the natural imbedding of a linear subshift,
 then
$(\bF[\Gamma])^r\cong
[\sum_\Gamma^0]^r\stackrel{i^*}{\rightarrow}\widehat{V}$
is a surjective $\bF[\Gamma]$-module homomorpism that is
$\widehat{V}$ is a finitely generated left $\bF[\Gamma]$-module.
Conversely, the dual of a finitely generated
$\bF[\Gamma]$-module is a linear subshift. It is important
to note that if $V\subset [\sum_\Gamma]^r$, \,$W\subset [\sum_\Gamma]^s$
are isomorphic linear subshifts then the dual of this isomorphism
provides a module-isomorphism between $\widehat{W}$ and $\widehat{V}$.
Conversely, the duals of isomorphic modules are isomorphic linear
subshifts.\qed
\section{The Noether property of group algebras}
Let $V\subset [\sum_\Gamma]^r$ be a linear subshift. We
 denote by $V^0$ the subspace of finitely supported elements.
\begin{proposition}
If $V^0$ contains a non-zero element then $\hg(V^0)>0$.
\end{proposition}
\proof
Let us suppose that a ball $B_r(\1g)$ contains the support
of a non-zero element in $V^0$. We claim that there exists
an $\epsilon>0$ such that if $n$
large enough then $\cf_n$ contains at least $\epsilon |\cf_n|$ disjoint
translates of $B_r(\1g)$. First note that the claim implies our
Proposition. If we have $M_n$ translates of $B_r(\1g)$ in $\cf_n$
then we can find $2^{M_n}$ different elements of $V^0$ which are all
supported in $\cf_n$. Therefore
$$\frac{\l2 |V_{\cf_n}|} {|\cf_n|}\geq \frac {M_n} {|\cf_n|}
\geq \frac {\epsilon |\cf_n|} {|\cf_n|}=\epsilon\quad.$$
Hence $\hg(V^0)\geq \epsilon$. Let us prove the claim.
Pick a maximal $2r$-net $a_1,a_2,\dots a_{M_n}$, that is 
maximal set of points in $\cf_n$ such that any two has distance 
greater or equal than $2r$. Then the $4r$-balls around
 the points $a_i$ are covering $\cf_n$. Hence
$M_n\geq  \frac{|\cf_n|} {B_{4r} (\1g)}$.
Then at least half of the $a_i$'s are not in $\cf_n\backslash
B_{r+1} (\partial\cf_n)$. The balls
around these elements far being from the boundary are completely
in $\cf_n$. Hence we have at least $\frac {1} {2} \frac{|\cf_n|} {B_{4r}
(\1g)}$ disjoint translates of $B_r(\1g)$ in $\cf_n$.\qed

\noindent
In the rest of this section we shall have an extra
assumption  on the amenable group $\Gamma$. We call an amenable
group Noether if $\bF[\Gamma]$ is a Noether ring. That is any
left submodule of $(\bF[\Gamma] )^r$ is finitely 
generated. According to Hall's theorem \cite{Pas} if $\Gamma$
is polycyclic-by-finite then $\Gamma$ is  Noether.
\begin{proposition}
\label{271}
If $V\subset [\sum_\Gamma]^r$ is a linear subshift and
$\Gamma$ is Noether, then $\hg(V)+\hg(V^\perp)=r$.
\end{proposition}
First of all $V^\perp\subset [\sum_\Gamma^0]^r\subset [\sum_\Gamma]^r$
hence the expression $\hg(V^\perp)$ is meaningful.
By our assumption $V^\perp$ is a finitely generated module, so let us
choose a $r_1,r_2\dots r_k$ finite generator set .
We need to prove that
$\lim_{n\rightarrow\infty}\frac {\l2 |V_{\cf_n}^\perp|} {|\cf_n|}=
r-\hg(V)$. In order to do so in it is enough to see that
\begin{equation}
\label{281}
\lim_{n\rightarrow\infty}\frac{\l2 |V_{\cf_n}^\perp|-\l2 |V_{n}^\perp|}
{|\cf_n|}=0\,,
\end{equation}
where $V_n^\perp$ denote the set of elements in $V^\perp$, supported in
$\cf_n$. Remember that $V_{\cf_n}^\perp$ denotes the restrictions of
 the elements of $V^\perp$, therefore $V_{\cf_n}^\perp\supset V_{n}^\perp$.
 By linear algebra,
$$\dim_{\bF}(V_{\cf_n})+\dim_{\bF}(V^\perp_{n})=r|\cf_n| $$
that is $\lim_{n\rightarrow\infty} \frac {\l2 |V_{n}^\perp|} {|\cf_n|}= 
r-\hg(V)$. Let us prove (\ref{281}). Any element of $V^\perp$ can be
written (not in a unique way !) in the form of $\sum_{i=1}^k a_ir_i$, where
$a_i\in\bF[\Gamma]$. Denote by $D$ the supremum of the diameters of
 the $r_i$'s.  If $supp (a_i)\subset\cf_n\backslash B_{D+1}(\partial\cf_n)$,
 for all $i$, then $\sum_{i=1}^k a_ir_i\in V_n^\perp$.
On the other hand if $supp(a_i)\cap B_{D+1}(\cf_n)=0$ for all $i$, then 
$\sum_{i=1}^k a_ir_i\mid_{\cf_n}=0$.
Therefore we have the pigeon-hole estimate
$$|V_{\cf_n}^\perp|\leq |V_n^\perp| 2^{kr B_{D+1} (\partial\cf_n)}$$
that is
$$\frac{\l2 |V_{\cf_n}^\perp|-\l2 |V_{n}^\perp|}
{|\cf_n|}\leq\frac{kr B_{D+1} (\partial\cf_n)}{|\cf_n|}$$
and the right hand side tends to zero.\qed

\noindent
Now we prove the density property.
\begin{proposition}
\label{291}
If $\Gamma$ is Noether and $V\subset[\sum_\Gamma]^r$ is
a linear subshift, then $\hg(V)=\hg(V^0)$.
\end{proposition}
\proof
By our previous Proposition, $\hg(V^\perp)=r-\hg(V)$. Therefore
$\hg(\overline{V^\perp})=r-\hg(V)$, where
$\overline{V^\perp}$ is the closure of $V^\perp$ as $[\sum_\Gamma^0]^r$
imbeds into $[\sum_\Gamma^0]^r$. Using our previous Proposition again,
$$\hg((\overline{V^\perp})^\perp)=\hg(V)\,.$$
 If $\xi\in (\overline{V^\perp})^\perp$ then $\xi$ is finitely
supported and $\xi\in(V^\perp)^\perp=V$, that is
$\xi\in V^0$. Therefore,
$$\hg(V)=\hg((\overline{V^\perp})^\perp)\leq\hg(V^0)\leq \hg(V)\quad\qed$$

\noindent
Actually our proofs of the last two Propositions gives a little bit
stronger result:
\begin{proposition}
\label{301}
Let $\Gamma$ be Noether and let $V$ be a linear subshift
such that $V^0$ is generated by $r_1,r_2,\dots,  r_k$ as left $\bF[\Gamma]$
-module. Denote by $\widetilde{V}_n$ the set of those elements in $V^0$
which can be written in the form of $\sum_{i=1}^k a_ir_i$, where
all the $a_i$'s are supported in $\cf_n$. Then  $\lim_{n\rightarrow\infty}
\frac{\l2 |\widetilde{V}_n|}{|\cf_n|}=\hg(V)$.
\end{proposition}

\section{The Yuzvinskii formula}
Recall Yuzvinskii's additivity formula for Abelian groups \cite{KS}.
Let $\Gamma\cong\bz^d$ and $\alpha$ be a $\Gamma$-action
 of continuous automorphisms on a compact metric group $X$. 
Suppose that $Y$ is a compact $\alpha$-invariant subgroup then
\begin{equation}
\label{381}
\hta(X)=\hta(Y)+\hta(X/Y)
\end{equation}
The results of Ward and Zhang \cite{WZ} suggest that a similar statement
might be true for general amenable actions. In our paper
we prove only a very special case.
\begin{proposition}
\label{382}
Let $Y\subset X\subset [\sum_\Gamma]^r$ be linear subshifts
where $\Gamma$ is Noether. Then
$$\hg(Y)=h_L^{top}(Y)+h_L^{top}(X/Y)=h_L^{top}(X)=\hg(X)\,,$$
where $L$ is the usual left $\Gamma$-action.
\end{proposition}
\proof
The key observation is the following lemma.
\begin{lemma}
\label{391}
Let $V\subset [\sum_\Gamma]^r$ be a linear subshift
where $\Gamma$ is Noether. Then there exists a constant $D$
such that if for some $\xi\in[\sum_\Gamma]^r$ with $\xi\mid_{B_D(\gamma)}\in
V_{B_D(\gamma)}$ for all $\gamma\in \Gamma$, then $\xi\in V$. That is
for Noether groups linear subshitfs are of finite type.
\end{lemma}
\proof
Let $V^\perp\subset[\sum_\Gamma^0]^r$ be the orthogonal
ideal of $V$. It is generated by $r_1,r_2,\dots ,r_N$, where
all $r_i$'s are supported in $B_D(\1g)$. then $\xi\notin V$
if and only if $\langle \xi,L_\gamma(r_i)\rangle\neq 1$
for some $i$ and $\gamma\in \Gamma$. It means that $\xi\mid_{B_D(\gamma)}
\notin V_{B_D(\gamma)}$.\qed

\noindent
Now we define a metric on $X/Y$ if $v,w\in X$ let $d([v],[w])=2^{-(n-1)}$,
where $n$ is the smallest integer such that
$(v-w)\mid_{B_D(\gamma)}\notin Y_{B_D(\gamma)}$ for all $\gamma\in\cf_n$.
here $D$ denotes the diameter of the joint support of a generator system
$s_1,s_2,\dots,s_M$ of the ideal $Y^\perp$.

\begin{lemma}
\label{401}
The metric  $d$ defines the pointwise convergence topology.
\end{lemma}

\proof
We need to prove that $d([v_n],0)\rightarrow 0$ implies
that $[v_n]\rightarrow Y$ in the factor topology of $X/Y$.
(the converse is obvious). Suppose that $\{[v_n]\}$ does not
converge to Y in the factor topology. Then there exists a subsequence
$v_{n_k}$ such that $v_{n_k}\rightarrow v\notin Y$
in the convergence topology of $[\sum_\Gamma]^r$. But then there exists a ball
$B_D(\gamma)$ such that $v_{n_k}\mid_{B_D(\gamma)}\notin Y_
{B_D(\gamma)}$ for large $k$. This contradicts to the
assumption that $d([v_{n_k}],0)\rightarrow 0$.\qed

\noindent
Now let us turn back to the proof of our Proposition. Similarly
to Proposition \ref{331} we have a lower estimate for $s_{X/Y}(n,1)$
in the $d$-metric. Let us denote by $G_n$ the set of elements in $Y^\perp$
which can be written in the form $\sum_{i=1}^m c_i s_i$ such that
 all the $c_i$'s are supported in $\cf_n$.
Denote by $H_n$ the set of those elements of $[\sum_\Gamma]^r$ which
are supported on $B_D(\cf_n)$ and orthogonal to $G_n$. Then by Propositions
\ref{291} and \ref{301}, $\frac {\l2 |G_n|}{|\cf_n|}\rightarrow \hg(Y)$.
We have the following inequality:
$$k_n= |X_{B_D(\cf_n)}|/|G_n\cap \,X_{B_D(\cf_n)}|\leq s_{X/Y}(n,1)\,.$$
Indeed, there exists $k_n$ elements of $X_{B_D(\cf_n)}$ such that their
pairwise differences $x_i-x_j\notin G_n$ thus
$\langle x_i-x_j,L_\gamma s_k\rangle\neq 1$ for some $s_k$ and $\gamma\in \cf_n$.
 Hence $d(L_{\gamma^{-1}}([x_i]),L_{\gamma^{-1}}([x_j])=1$.
Consequently, $\hg(X)-\hg(Y)\leq h^{top}_L(X/Y)$. Now fix an $\epsilon$ and
 let $B_R(\1g)\supset\cf_m$, where $2^{-m}<\epsilon$.
Then obviously,
$$s_{X/Y}(n,\epsilon)\leq\frac {|X_{B_{D+r}(\cf_n)}|} {|Y_{B_{D+r}(\cf_n)}|}\,,$$
which implies the converse inequality : $\hg(X)-\hg(Y)\geq h^{top}_L(X/Y)$.
\section{Betti numbers}
In this section we define an analogue of the $L^2$-Betti numbers. 
Let $\wk$ be a regular, normal $\Gamma$-covering of
a finite simplicial complex $K$, where $\Gamma$ is an amenable
group that acts freely and simplicially on $\wk$ and $\wk/\Gamma=K$.
We have the ordinary cochain complex of
$\bF$-coefficients over $\wk$:
$$ C^0(\wk,\bF)\stackrel{d_0}{\rightarrow}
C^1(\wk,\bF)\stackrel{d_1}{\rightarrow}\dots\stackrel{d_{n-1}}{\rightarrow}
C^n(\wk,\bF)\,,$$
Then the $p$-cochain space $C^p(\wk,\bF)$ is $\Gamma$-isomorphic
to $[\sum_\Gamma]^{|K_p|}$, where $K_p$ denotes
the set of $p$-simplices in $K$.
We define the $p$-th entropy Betti number $\ben^p(K)$ as
$\hg(Ker\,d_p)-\hg(Im\,d_{p-1})$. The following theorem is the
analogue of Cohen's theorem \cite{Coh}
\begin{proposition}
\label{431}
$\sum^n_{p=0} (-1)^p \ben^p(K)$ equals to the Euler-characteristics
of $K$.
\end{proposition}
\proof
By Proposition \ref{181}, 
$$b_E^p(K)=\hg(Ker\,d_p)+\hg(Ker\,d_{p-1})-\hg(C^{p-1}(\wk,\bF))\,.$$
Summing up these equations for all $p$ with alternating signs we obtain
that
$$\sum^n_{p=0} (-1)^p \ben^p(K)=\sum^n_{p=0} (-1)^p |K_p|=e(K)\,\qed$$

\noindent
Now let us see the analogue of the result of Cheeger \& Gromov.
\begin{proposition}
\label{432}
If $\wk$ is contractible, then all entropy Betti numbers are vanishing.
\end{proposition}
The proof is much easier than for the $L^2$-Betti numbers.
 If $p>0$, then $\bF$-cohomologies are vanishing 
therefore $\ben^p(K)=0$ if $p>0$. If $p=0$, then the cocycle space is finite
so the entropy Betti number must be zero. \qed

\noindent
\begin{corollary}
If $K$ is a finite acyclic simplicial complex with an amenable
 fundamental group then its Euler characteristics is zero.
\end{corollary}
Now we prove the analogue of the result of Dodziuk and Mathai.
\begin{proposition}\label{najo}
$$\ben^p(K)=\lim_{n\rightarrow\infty}\frac
{\dim_{\bF} H^p(L_n, \bF)}{|\cf_n|}\,,$$
where $\{L_n\}$ is an exhaustion of $\wk$ spanned by the \Fol-sets.
\end{proposition}
{\bf Remark:} Since $\dim_{\bF} H^p(L_n, \bF)\geq \dim_{\br} H^p(L_n, \br)$
the entropy Betti numbers are at least as large as the corresponding
$L^2$-Betti numbers. It is easy to construct examples, where some entropy
Betti numbers are strictly larger than the corresponding $L^2$-Betti number
(cf. the remark after Proposition \ref {482}.)

\noindent
{\bf Proof:} (of Proposition \ref{najo}) 
First note again that $C^0(\wk,\bF)\cong[\sum_\Gamma]^{|K_p|}$. Denote by
 $R$ a constant such that for any $[(\gamma,p),(\delta,q)]$ 1-simplex
of $\wk$, $d_\Gamma(\gamma,\delta)\leq R$. Now we can build an
extended configuration for $Ker\,d_p$ and $Im\,d_p$ the
following way. Let $S(\Lambda)$ be the simplicial complex spanned
by vertices of the form $(\gamma,p)$, where $\gamma\in B_{2r}(\Lambda)$
and $p\in K_0$. Also, let $L_n$ be the simplicial complex spanned
by the vertices with first coordinate in $\cf_n$. Consider the coboundary
operator as $[\sum_\Gamma]^{|K_p|}\stackrel{d_p}{\rightarrow} 
[\sum_\Gamma]^{|K_{p+1}|}$. Let $A_p(\Lambda)$ be the space of those
functions in $[\sum_\Gamma]^{|K_p|}$ which are supported on $\Lambda$
and are the restriction of a cocycle of $S(\Lambda)$
respectively let $B_p(\Lambda)$ be
the space of restrictions of coboundaries of $S(\Lambda)$.
 Obviously,
$A_p(\Lambda)$ is an extended configuration of $Ker\,d_p$ and
$B_p(\Lambda)$ is an extended configuration of $Im\,d_p$. Then the
usual pigeon-hole argument and Proposition \ref {151}. implies that
$$\hg(Ker\,d_p)=\lim_{n\rightarrow\infty}\frac
{ \dim_{\bF} (Z^p(S(\cf_n)))}{|\cf_n|}$$
$$\hg(Im\,d_p)=\lim_{n\rightarrow\infty}\frac
{ \dim_{\bF} (B^p(S(\cf_n)))}{|\cf_n|}\,,$$
where $Z^p$ resp. $B^p$ denote the space of cocycles resp. coboundaries.
Therefore
$$\ben^p(K)=\lim_{n\rightarrow\infty}\frac
{ \dim_{\bF} (H^p(S(\cf_n),\bF))}{|\cf_n|}\,.$$
Finally we must prove that
$$\lim_{n\rightarrow\infty}\frac
{ \dim_{\bF} (H^p(S(\cf_n),\bF))-\dim_{\bF}(H^p(L_n,\bF))}{|\cf_n|}=0\,.$$
Note that it follows from the long exact cohomology sequence induced by
the inclusion $L_n\rightarrow S(\cf_n)$ and the obvious fact that
$\frac{\dim_{\bF} (H^p(S(\cf_n),L_n,\bF)}{|\cf_n|}$
tends to zero as $n\rightarrow\infty$.\qed
\section{Towers and fixed points}
In this section we recall some ideas of Farber \cite{Far}. Let $\Gamma$
be a finitely generated residually-$p$ group. That is, there exists
a chain of normal subgroups of prime power index,
$\Gamma\supset\Gamma_1\supset\Gamma_2\supset\dots\quad,\mbox{where}\,
\cap_{j=1}^\infty\Gamma_j=\{1_\Gamma\}$.
Let $\wk/\Gamma=K$ be as in the previous section. Then
one can consider the tower of finite simplicial complexes
$X_i=\wk/\Gamma_i$. Note that
$X_i$ is a simplicial $(\Gamma:\Gamma_i)$-covering
of $K$. Farber proved (\cite{Far}, Theorem 11.1 that
$\lim_{j\rightarrow\infty}\frac
{\dim_{\bF} H^i(X_j,\bF)}{|\Gamma:\Gamma_i|}$ always
exists. The
following conjecture is the  analogue of L\"uck's
theorem on approximating the $L^2$-Betti numbers \cite{Luck2}:
\begin{conjecture}
\label{481}
If $\Gamma$ is as above a residually-$2$ group, then
$$\lim_{j\rightarrow\infty}\frac
{\dim_{\bF} H^i(X_j,\bF)}{|\Gamma:\Gamma_i|}=\ben^i(K)\,.$$
\end{conjecture}
\begin{proposition}
\label{482}
The conjecture is true if $\Gamma$ is free Abelian.
\end{proposition}
\proof
First of all note that if $\Gamma$ is Noether, then any
$V\subset[\sum_\Gamma]^r$
linear subshift is expansive. That is there exists $\epsilon>0$ such
 that if $x\neq y\in V$, then for some
$\gamma\in\Gamma$ : $d(L_\gamma(x),L_\gamma(y))\geq\epsilon$. This is just
a reformulation of Lemma \ref{391}. The following result is due to
Schmidt (\cite{KS}, Theorem 21.1).
\begin{proposition}
If $\alpha$ is an expansive $\bz^d$-action by automorphisms of a compact
Abelian group $X$, then
$$\lim_{|\bz^d:\Lambda|\rightarrow\infty,\,\Lambda\subset\bz^d}
\frac{|Fix\,\Lambda|}{|\bz^d:\Lambda|}=\hta(X)\,$$
where $Fix\,\Lambda$ denotes the set of fixed points of the subgroup
$\Lambda$.
\end{proposition}
Now let us turn to the proof of Proposition \ref{482}.
Let $Z^i_j$ be the space of $i$-cocycles on $X_j$ and $Z^i$ be the
space of $i$-cocycles on $\wk$. Then $Z^i_j$ is exactly the set of
fixed points of the subgroup $\Gamma_j$ on $Z^i$. (Note that
the similar statement on coboundaries would not be necessarily true).
By Proposition \ref{481}, $\lim_{j\rightarrow\infty}\frac
{\dim_{\bF} Z^i_j}{|\Gamma:\Gamma_j|}=\hg(Z^i)$.
If $C^i_j$ denotes the space of $i$-cochains on $X_j$
and $C^i$ denotes the space of $i$-cochains on $\wk$, then
$\frac{\dim_{\bF} C^i_j}{|\Gamma:\Gamma_j|}=|K_i|=\hg(C^i)$,
for all $j$. By our Proposition \ref{181},
$$\ben^i(K)=\hg(Z^i)+\hg(Z^{i-1})-\hg(C^{i-1})\,.$$
Also, $$\dim_{\bF} H^i(X_j,\bF)=\dim_{\bF}(Z^i_j)+\dim_{\bF}(Z^{i-1}_j)-
\dim_{\bF}(C^{i-1}_j)\,.$$  Hence our Proposition follows.\qed

\noindent
It is not hard to construct a $\wk$, where for
some $p$ the entropy and $L^2$-Betti numbers differ.
Simply consider the Cayley graph of $\bz^d$ and
then just stick a $\br P^4$ on each vertex.
Then if $L_n$ denote the approximative complexes for
 some \Fol-exhaustion :
$$b_E^4(K)=\lim_{n\rightarrow\infty}\dim_{\bF}\frac{H^4(L_n,\bF)}{|\cf_n|}=1$$
and
$$\bet^4(K)=\lim_{n\rightarrow\infty}\dim_{\br}\frac{H^4(L_n,\br)}
{|\cf_n|}=0\,.$$
\section{The Grothendieck group and the integrality of the Betti numbers}
First we recall
the notion of the Grothendieck group of a non-commutative ring $R$
\cite{Ros}.
Let $G(R)$ be the Abelian group, defined by generators
$\{[M]\}$, where the $M$'s are the finitely generated 
left $R$-modules up to isomorphism. The relations
are in the form \,$[M]+[N]=[L]$, for any exact
sequence $0\rightarrow M\rightarrow L\rightarrow N\rightarrow 0$.
L\"uck \cite{Luck} proved that if $R=\bc[\Gamma]$, where $\Gamma$
is amenable, then $[\bc[\Gamma]]$ generates an infinite cyclic subgroup in
$G(R)$.
\begin{proposition}
$[\bc[\Gamma]]$ generates an infinite cyclic subgroup
in $G(\bc[\Gamma])$ for any finitely generated amenable group $\Gamma$.
\end{proposition}
\proof
It is enough to define a rank on finitely generated $\bc[\Gamma]$-modules,
such that $rk([\bc[\Gamma]])=1$ and
$rk([M])+rk([N])=rk([L])$ if
$$0\rightarrow M\stackrel{i}{\rightarrow} L \stackrel{p}{\rightarrow}N
\rightarrow 0\,.$$
Let $rk(M)=\hg(\widehat{M})$. Now apply Proposition \ref{181} for the
subshifts
$$0\rightarrow \widehat{N}\stackrel{\widehat{p}}{\rightarrow}
\widehat{L}\stackrel{\widehat{i}}{\rightarrow}
\widehat{M}
\rightarrow 0$$
and the additivity follows.\qed

\noindent
Linnell \cite{Lin} proved that all $L^2$-Betti
numbers are integers for torsion-free elementary amenable group $\Gamma$.
We can prove the following proposition.
\begin{proposition}
If $\Gamma$ is poly-infinite-cyclic, then $\hg(V)$ is an integer
for any linear subshift $V$.
\end{proposition}
\proof
Let $M=\widehat{V}$ be the dual $\cf_2[\Gamma]$-module
of our subshift. Then, by Theorem 3.13 \cite{Pas} $M$ has a
finite resolution by finitely generated projective modules :
$$0\rightarrow M_n\rightarrow\dots\rightarrow M_2\rightarrow M_1\rightarrow
 M\rightarrow  0$$
Then, as we pointed out earlier the dual sequence
$$0\rightarrow V\rightarrow V_1\rightarrow V_2\rightarrow
\dots \rightarrow V_n\rightarrow 0$$ 
is an exact sequence of linear subshifts and continuous homomorphisms,
where $V_i=\widehat{M_i}$. By Proposition \ref{181} it is enough to show
that all the $\hg(V_i)$'s are integers. By a result of Grothendieck \&
Serre (Theorem 4.13 \cite{Pas}) if $\Gamma$ is poly-infinite-cyclic,
then all finitely generated, projective $\bF[\Gamma]$-module is stably
free. Hence, using the notation of the previous section :
$$\hg(V_i)=rk(\widehat{V_i})=rk((\bF[\Gamma])^n)-rk((\bF[\Gamma])^m)=n-m$$
is an integer. \qed

\end{document}